\documentclass[11pt]{article}
\usepackage{amsmath,amssymb,amsthm,mathrsfs,graphicx,color,epsfig,verbatim}

\newtheorem{thm}{Theorem}[section]

\newtheorem{defn}{Definition}
\newtheorem{prop}[thm]{Proposition}
\newtheorem{lemma}[thm]{Lemma}
\newtheorem{cor}[thm]{Corollary}

\newcommand{\mpf}{\mathcal{MP}_G}
\newcommand{\fG}{\mathcal{F}_G}
\newcommand{\gmp}{$G$-multiparking function} 
\newcommand{\gmps}{$G$-multiparking functions}

\begin{document}

\title{Multiparking Functions, Graph Searching, and the Tutte Polynomial}

\author{
Dimitrije Kosti\'{c}$^{1}$\thanks{Partially supported by Texas A\&M's NSF
VIGRE grant.} \: and \: Catherine H. Yan$^{2}$\thanks{Partially support
  by NSF grant DMS 0245526.} 
\vspace{.3cm} \\
$^{1,2}$ Department of Mathematics \\
Texas A\&M University,  College Station, TX 77843\\
$^{2}$ Center for Combinatorics, LPMC\\
Nankai University, Tianjin 300071,  P.R. China 
\vspace{.4cm} \\
$^1$ dkostic@math.tamu.edu,   $^2$cyan@math.tamu.edu
} 
\date{} 
\maketitle



\begin{abstract}

A parking function of length $n$ is a sequence $(b_1, b_2, \dots,
 b_n)$ of nonnegative integers for which there is a permutation $\pi \in S_n$ so that $0 \leq b_{\pi(i)} < i$ for all $i$.  A well-known result 
about parking functions is that the polynomial $P_n(q)$, which enumerates the complements of parking functions by the sum of their terms, is the 
generating function for the number of connected graphs by the number of excess edges when evaluated at $1+q$. In this paper we extend this result to 
arbitrary connected graphs $G$. In general the polynomial that encodes information about subgraphs of $G$ is the Tutte polynomial $t_G(x,y)$, which 
is the generating function for two parameters, namely the internal and external activities, associated with the spanning trees of $G$.  We define 
\gmps, which generalize the $G$-parking functions that Postnikov and Shapiro introduced in the study of certain quotients of the polynomial ring.  
We construct a family of algorithmic bijections between the spanning forests of a graph $G$ and the $G$-multiparking functions.  In particular, the 
bijection induced by the breadth-first search leads to a new characterization of external activity, and hence a representation of Tutte polynomial 
by the reversed sum of \gmps.

\end{abstract}

\section{Introduction}

The (classical) parking functions of length $n$ are sequences $(b_1, b_2, \ldots, b_n)$ of nonnegative integers for which there is a permutation 
$\pi \in S_n$ so that $0 \leq b_{\pi(i)} < i$ for all $i$.  This notion was first introduced by Konheim and Weiss \cite{KonheimWeiss:1966} in the 
study of the linear probes of random hashing function. The name comes from a picturesque description in \cite{KonheimWeiss:1966} of the sequence of 
preferences of $n$ drivers under certain parking rules.  Parking functions have many interesting combinatorial properties.  The most notable one is 
that the number of parking functions of length $n$ is $(n+1)^{n-1}$, Cayley's formula for the number of labeled trees on $n+1$ vertices. This 
relation motivated much work in the early study of parking functions, in particular, combinatorial bijections between the set of parking functions 
of length $n$ and labeled trees on $n+1$ vertices.  See \cite{GiKa:1999} for an extensive list of references.

There are a number of generalizations of parking functions, for example, see \cite{CoriPoulalhon:2002} for the double parking functions, 
\cite{Stanley:1998,Yan:1997,Yan:2001} for $k$-parking functions, and \cite{PitmanStanley:2002,KungYan:2003} for parking functions associated with an 
arbitrary vector. Recently, Postnikov and Shapiro \cite{PostnikovShapiro:2004} proposed a new generalization, the {\em $G$-parking functions}, 
associated to a general connected digraph $D$. Let $G$ be a digraph on $n+1$ vertices indexed by integers from $0$ to $n$.  A $G$-parking function 
is a function $f$ from $[n]$ to $\mathbb N$, the set of non-negative integers, satisfying the following condition: for each subset $U \subseteq [n]$ 
of vertices of $G$, there exists a vertex $j \in U$ such that the number of edges from $j$ to vertices outside $U$ is greater than $f(j)$.  For the 
complete graph $G=K_{n+1}$, such defined functions are exactly the classical parking functions, where one views $K_{n+1}$ as the digraph with one 
directed edge $(i,j)$ for each pair $i \neq j$.  In \cite{ChPy:2005} Chebikin and Pylyavskyy constructed a family of bijections between the set of 
$G$-parking functions and the (oriented) spanning trees of that graph.

Perhaps the most important statistic of the classical parking functions is the (reversed) sum, that is, ${n \choose 2}-(x_1+x_2+\cdots +x_n)$ for a 
parking function $(x_1, x_2,\dots, x_n)$ of length $n$. It corresponds to the number of linear probes in hashing functions \cite{Knuth98}, the 
number of inversions in labeled trees on $[n+1]$ \cite{Kreweras}, and the number of hyperplanes separating a given region from the base region in 
the extended Shi arrangements \cite{Stanley:1998}, to list a few.  It is also closely related to the number of connected graphs on $[n+1]$ with a 
fixed number of edges. In \cite{Yan:2001} the second author gave a combinatorial explanation, which revealed the underlying correspondence between 
the classical parking functions and labeled, connected graphs. The main idea is to use breadth-first search to find a labeled tree on any given 
connected graph, and record such a search by a queue process.

The objective of the present paper is to extend the result of \cite{Yan:2001} to arbitrary graphs. For a general graph $G$, a suitable tool to study 
all subgraphs of $G$ is the Tutte polynomial.  This is a generating function with two parameters, the internal and external activities, which are 
functions on the collection of spanning trees of $G$.  Evaluating the Tutte polynomial at various points provides combinatorial information about 
the graph, for example, the number of spanning trees, spanning forests, connected subgraphs, acyclic orientations, subdigraphs, etc.  These many 
valuations make the Tutte polynomial one of the most fundamental tools in algebraic graph theory.
 
An important approach to getting information about the Tutte polynomial is to use partitions.  This approach dates from the 1960's, see Crapo 
\cite{Crapo:1969}. More information on the history of Tutte polynomial can be found in \cite{GesselSagan:1996}.  Also in that paper Gessel and Sagan 
proposed a number of new notions of external activity, along with a new way to partition the substructures of a given graph. The basic method is to 
use depth-first search to associate a spanning forest $F$ with each substructure to be counted.  This process partitions the simplicial complex of 
all substructures (ordered by inclusion) into intervals, one for each $F$. Every interval turns out to be a Boolean algebra consisting of all ways 
to add external active edges to $F$. Expressing the Tutte polynomial in terms of sums over such intervals permits one to extract the necessary 
combinatorial information.

In \cite{GesselSagan:1996} Gessel and Sagan also mentioned another search, the {\em neighbors-first search}, and related the external activity 
determined by the neighbors-first search on a complete graph with $n+1$ vertices to the sum of (classical) parking functions of length $n$.  This 
connection was further explained in \cite{Yan:2001}.  In the present paper we extend this result to an arbitrary graph $G$ by developing the 
connection between Tutte polynomial of $G$ and certain restricted functions defined on $V(G)$, the vertex set of $G$. This is achieved by combining 
the two approaches mentioned before.  First, we use breadth-first search to get a new partition of all spanning subgraphs of $G$. Each subgraph is 
associated with a spanning forest of $G$, which allows us to get a new expression of the Tutte polynomial in terms of breadth-first external 
activities of its spanning forests.  Second, we construct bijections between the set of all spanning forests of $G$ and the set of functions defined 
on $V(G)$ with certain restrictions. One of such bijection, namely the one induced by breadth-first search with a queue, leads to the 
characterization of the (breadth-first) external activity of a spanning forest by the corresponding function.

To work with spanning forests, we propose the notion of a $G$-multiparking function, a natural extension of the notion of a $G$-parking function. 
Let $G$ denote a graph with a totally ordered vertex set $V(G)$.  Often we will take $V(G)=[n]=\{1, 2, \dots, n\}$.  For simplicity and clarity, we 
assume that $G$ is a simple graph in most of the paper, except at the end of Section 3 where we explain how our construction could be modified to 
apply for general directed graphs, with possible loops and multiple edges. This includes undirected graphs as special cases, as an undirected graph 
can be viewed as a digraph where each edge $\{u,v\}$ is replaced by a pair of arcs $(u,v)$ and $(v,u)$.

For any subset $U \subseteq V(G)$, and vertex $v \in U$, define $outdeg_U(v)$ to be the cardinality of the set $\{ \{v,w\} \in E(G) | w \notin U 
\}$. Here $E(G)$ is the set of edges of $G$.

\begin{defn} Let $G$ be a simple graph with $V(G)=[n]$.  A \emph{G-multiparking function} is a function $f : V(G)= [n]\rightarrow \mathbb{N} \cup \{ 
\infty \}$, such that for every $U \subseteq V(G)$ either \textrm{\bf (A)} $i$ is the vertex of smallest index in $U$, (written as $i=\min(U)$), and 
$f(i) = \infty$, or \textrm{\bf (B)} there exists a vertex $i \in U$ such that $0 \leq f(v_i) < outdeg_U(i)$. \end{defn}

The vertices which satisfy $f(i)= \infty$ in \pmb{(A)} will be called \emph{roots of $f$} and those that satisfy \pmb{(B)} (in $U$) are said to be 
\emph{well-behaved} in $U$, and \pmb{(A)} and \pmb{(B)} will be used to refer, respectively, to these conditions hereafter.  Note that vertex $1$ is 
always a root. The $G$-multiparking functions with only one root (which is necessarily vertex $1$) are exactly the $G$-parking functions, as defined 
by Postnikov and Shapiro.

Sections 2 and 3 are devoted to the combinatorial properties of $G$-multiparking functions.  In \S2 we construct a family of algorithmic bijections 
between the set $\mathcal{MP}_G$ of $G$-multiparking functions and the set $\mathcal{F}_G$ of spanning forests of $G$.  Each bijection is a process 
based on a choice function, (c.f. \S2), which determines how the algorithm proceeds.  In \S3 we give a number of examples to illustrate various 
forms of the bijection.  This includes the cases where there is a special order on $V(G)$, for instance {\em depth-first search order}, {\em 
breadth-first search order}, and a prefixed linear order, and the cases that the process possesses certain data structure, such as queue and stacks.  
At the end of \S3 we explain how the algorithm works for general directed graphs.

Section 4 is on the relation between $G$-multiparking functions and the Tutte polynomial of $G$. First, for each forest $F$ we give a 
characterization of $F$-redundant edges, which are edges of $G-F$ that are ``irrelevant'' in determining the corresponding \gmp. Using that we 
classify the edges of $G$, and establish an equation between $|E(G)|$ and sum of \gmp, $|E(F)|$, and the $F$-redundant edges.  Then we use 
breadth-first search to partition all the subgraphs of $G$ into intervals.  Each interval consists of all graphs obtained by adding some 
breadth-first externally active edges to a spanning forest $F$. The set of breadth-first externally active edges of $F$ are exactly the 
$F$-redundant edges of a certain type, which allows us to express the number of breadth-first externally active edges, and hence the Tutte 
polynomial, by the values of corresponding \gmp.  In section 5 we exhibit some enumerative results related to \gmps \ and substructures of graphs.

\section{Bijections between multiparking functions and spanning forests}

In this section, we construct bijections between the set $\mathcal{MP}_G$ of $G$-multiparking functions and the set $\mathcal{F}_G$ of spanning 
forests of $G$.  For simplicity, here we assume $G$ is a simple graph with $V(G)=[n]$.  A \emph{sub-forest} $F$ of $G$ is a subgraph of $G$ without 
cycles.  A leaf of $F$ is a vertex $v \in V(F)$ with degree 1 in $F$.  Denote the set of leaves of $F$ by $Leaf(F)$.  Let $\prod$ be the set of all 
ordered pairs $(F, W)$ such that $F$ is a sub-forest of $G$, and $\emptyset \neq W \subseteq Leaf(F)$. A {\em choice function} $\gamma$ is a 
function from $\prod$ to $V(G)$ such that $\gamma(F, W) \in W$.  Examples of various choice functions will be given in \S3, where we also explain 
how the bijections work on a general directed graph, in which loops and multiple edges are allowed. As one can see, loopless undirected graphs can 
be viewed as special case there.

Fix a choice function $\gamma$. Given a $G$-multiparking function $f \in \mathcal{MP}_G$, we define an algorithm to find a spanning forest $F \in 
\mathcal{F}_G$.  Explicitly, we define quadruples $(val_i, P_i, Q_i, F_i)$ recursively for $i=0, 1, \dots, n$, where $val_i: V(G) \rightarrow 
\mathbb{Z}$ is the {\em value function}, $P_i$ is the set of \emph{processed} vertices, $Q_i$ is the set of vertices {\em to be processed}, and 
$F_i$ is a subforest of $G$ with $V(F_i)=P_i \cup Q_i$, $Q_i \subseteq Leaf(F_i)$ or $Q_i$ consists of an isolated vertex of $F_i$.

\noindent {\bf Algorithm A}. 

\begin{itemize}

\item {\bf Step 1: initial condition.} Let $val_0 = f$, $P_0$ be empty, and $F_0=Q_0=\{1\}$.

\item {\bf Step 2: choose a new vertex $v$.} At time $i \geq 1$, let $v=\gamma(F_{i-1}, Q_{i-1})$, where $\gamma$ is the choice function.

\item {\bf Step 3: process vertex $v$.} For every vertex $w$ adjacent to $v$ and $w \notin P_{i-1}$, set $val_i(w)=val_{i-1}(w)-1$. For any other 
vertex $u$, set $val_i(u)=val_{i-1}(u)$.  Let $N=\{w | val_{i}(w)=-1, val_{i-1}(w) \neq -1 \}$.  Update $P_i$, $Q_i$ and $F_i$ by letting 
$P_i=P_{i-1} \cup \{v\}$, $Q_i=Q_{i-1} \cup N \setminus \{ v\}$ if $Q_{i-1} \cup N \setminus \{ v\} \neq \emptyset$, otherwise $Q_i=\{u\}$ where $u$ 
is the vertex of the lowest-index in $[n] - P_{i}$. Let $F_i$ be a graph on $P_i \cup Q_i$ whose edges are obtained from those of $F_{i-1}$ by 
joining edges $\{w, v \}$ for each $w \in N$.  We say that the vertex $v$ is processed at time $i$.

\end{itemize}

Iterate steps 2-3 until $i=n$. We must have $P_n=[n]$ and $Q_n = \emptyset$.  Define $\Phi=\Phi_{\gamma, G}: \mathcal{MP}_G \rightarrow 
\mathcal{F}_G$ by letting $\Phi(f)=F_n$.

If an edge $\{v,w\}$ is added to the forest $F_i$ as described in Step 3, we say that {\em $w$ is found by $v$}, and $v$ is the \emph{parent} of 
$w$, if $v \in P_{i-1}$.  (In this paper, the parent of vertex $v$ will be frequently denoted $v^p$.)  By Step 3, a vertex $w$ is in $Q_i$ because 
either it is found by some $v$ that has been processed, and $\{v,w\}$ is the only edge of $F_i$ that has $w$ as an endpoint, or $w$ is the 
lowest-index vertex in $[n] - P_i$ and is an isolated vertex of $F_i$.  Also, it is clear that each $F_i$ is a forest, since every edge $\{u,w\}$ in 
$F_i \setminus F_{i-1}$ has one endpoint in $V(F_i)\setminus V(F_{i-1})$. Hence $\gamma(F_i, Q_i)$ is well-defined and thus we have a well-defined 
map $\Phi$ from $\mpf$ to $\fG$.  The following proposition describes the role played by the roots of a $G$-multiparking function $f$.

\begin{prop} \label{root} Let $f$ be a \gmp.  Each tree component $T$ of $\Phi(f)$ has exactly one vertex $v$ with $f(v) = \infty$. In particular, 
$v$ is the least vertex of $T$.
\end{prop} 

\begin{proof} In the algorithm A the value for a root of $f$ never changes, as $\infty -1 =\infty$. Each nonroot vertex $w$ of $T$ is found by some 
other vertex $v$, and $\{v, w\}$ is an edge of $T$.  As any tree has one more vertex than its number of edges, it has exactly one vertex without a 
parent. By the definition of Algorithm A, this must be a root of $f$.

To show that the root is the least vertex in each component, let $r_1< r_2< \dots < r_k$ be the roots of $f$ and suppose $T_{1}, T_{2}, \ldots, 
T_{k}$ are the trees of $F=\Phi(f)$, where $r_i \in T_i$.  Let $T_j$ be the tree of smallest index $j$ such that there is a $v \in T_j$ with $v < 
r_j$.  Then $j > 1$ since the vertex $1$ is always a root.  Define $U := V( T_{j} \cup T_{j+1} \cup \ldots \cup T_{k})$.  $U$ is thus a proper 
subset of $V(G)=[n]$.  By assumption, the vertex of least index in $U$ is not a root.  Therefore, $U$ must contain a well-behaved vertex; that is, a 
vertex $v$ such that $0 \leq f(v) < outdeg_U(v)$.  Note that all the edges counted by $outdeg_U(v)$ lead to vertices in the trees $T_{1}, T_{2}, 
\ldots, T_{j-1}$.  By the structure of algorithm $A$, all the vertices in the first $j-1$ trees are processed before the parent of $v$ is processed.  
But this means that by the time $A$ processes all the vertices in the first $j-1$ trees, $val_i(v) = f(v) - outdeg_U(v) \leq -1$, so $v$ should be 
adjacent to some vertex in one of the first $j-1$ trees.  This is a contradiction.

\end{proof} 

From the above proof we also see that the forest $F=\Phi(f)$ is built tree by tree by the algorithm A. That is, if $T_i$ and $T_j$ are tree 
components of $F$ with roots $r_i$, $r_j$ and $r_i < r_j$, then every vertex of tree $T_i$ is processed before any vertex of $T_j$.

To show that $\Phi$ is a bijection, we define a new algorithm to find a \gmp \ for any given spanning forest, and prove that it gives the inverse 
map of $\Phi$.

Let $G$ be a graph on $[n]$ with a spanning forest $F$. Let $T_1, \dots, T_k$ be the trees of $F$ with respective minimal vertices $r_1=1 < r_2 < 
\cdots < r_k$.

\noindent{\bf Algorithm B.} 

\begin{itemize}

\item {\bf Step 1. Determine the process order $\pi$.} Define a permutation $\pi=(\pi(1), \pi(2), \dots, \pi(n))= (v_1 v_2 \dots v_n)$ on the 
vertices of $G$ as follows.  First, $v_1=1$.  Assuming $v_1, v_2, \dots, v_i$ are determined,

\begin{itemize}
   
\item {\bf Case (1)} If there is no edge of $F$ connecting vertices in $V_i=\{v_1, v_2, \dots, v_i\}$ to vertices outside $V_i$, let $v_{i+1}$ be 
the vertex of smallest index not already in $V_i$;

\item {\bf Case (2)} Otherwise, let $W=\{v \notin V_i: v \text{ is adjacent to some vertices in } V_i\}$, and $F'$ be the forest obtained by 
restricting $F$ to $V_i \cup W$. Let $v_{i+1}=\gamma(F', W)$.

\end{itemize}

(Hereafter, when discussing process orders, we will write $v_i$ as $\pi(i)$.)

\item {\bf Step 2. Define a \gmp \ $f=f_F$.} Set $f(r_1)=f(r_2)=\cdots =f(r_k)=\infty$.  For any other vertex $v$, let $r_v$ be the minimal vertex 
in the tree containing $v$, and $v, v^p, u_1, \dots, u_t, r_v$ be the unique path from $v$ to $r_v$.  Set $f(v)$ to be the cardinality of the set 
$\{ v_j | ( v, v_j ) \in E(G), \; \pi^{-1}(v_j) < \pi^{-1}(v^p) \}$.

\end{itemize}

To verify that a function $f=f_F$ defined in this way is a \gmp, we need the following lemma.

\begin{lemma}\label{lemma_2.1} Let $f: V(G) \rightarrow \mathbb{N} \cup \{ \infty \}$ be a function.  If $v \in U \subseteq V(G)$ obeys property 
\textrm{\bf (A)} or property \textrm{\bf (B)} and $W$ is a subset of $U$ containing $v$, then $v$ obeys the same property in $W$.

\end{lemma}

\begin{proof} If $f(v)=\infty$ and $v$ is the smallest vertex in $U$, then clearly it will still be the smallest vertex in $W$.  If $v$ is 
well-behaved in $U$, then $0 \leq f(v) < outdeg_U(v)$ and as $W \subseteq U$, we have $outdeg_U(v) \leq outdeg_W(v)$.  Thus $v$ is well-behaved in 
$W$.

\end{proof}

The \emph{burning algorithm} was developed by Dhar \cite{Dhar:1990} to determine if a function on the vertex set of a graph had a property called 
\emph{recurrence}. An equivalent description for $G$-parking functions is given in \cite{ChPy:2005}: We mark vertices of $G$ starting with the root 
$1$. At each iteration of the algorithm, we mark all vertices $v$ that have more marked neighbors than the value of the function at $v$. The 
function is a $G$-parking function if and only if all vertices are marked when this process terminates.  Here we extend the burning algorithm to 
$G$-multiparking functions, and write it in a linear form.

\begin{prop}~\label{burning} A vertex function is a $G$-multiparking function if and only if there exists an ordering $\pi(1), \pi(2), \ldots, 
\pi(n)$ of the vertices of a graph $G$ such that for every $j$, $\pi(j)$ satisfies either condition \textrm{\bf (A)} or condition \textrm{\bf(B)} in 
$U_j := \{ \pi(j), \ldots, \pi(n)\}$.

\end{prop}

\begin{proof} 

We say that the vertices can be ``thrown out" in the order $\pi(1), \pi(2), \ldots, \pi(n)$ if they satisfy the condition described in the 
proposition. By the definition of \gmp, it is clear that for a \gmp, vertices can be thrown out in some order.

Conversely, suppose that for a vertex function $f: V(G) \rightarrow \mathbb{N} \cup \{ \infty \}$ the vertices of $G$ can be thrown out in a 
particular order $\pi(1), \pi(2), \ldots, \pi(n) $. For any subset $U$ of $V(G)$, let $k$ be the maximal index such that $U \subseteq U_k=\{ \pi(k), 
\dots, \pi(n)\}$. This implies $\pi(k) \in U$.  But $\pi(k)$ satisfies either condition \pmb{(A)} or condition \pmb{(B)} in $U_k$.  By Lemma 
\ref{lemma_2.1}, $\pi(k)$ satisfies either condition \pmb{(A)} or condition \textrm{\bf (B)} in $U$. Since $U$ is arbitrary, $f$ is a \gmp. 

\end{proof}

\begin{prop} The Algorithm B, when applied to a spanning forest of $G$, yields a \gmp \ $f=f_F$.

\end{prop} 

\begin{proof} Let $\pi$ be the permutation defined in Step 1 of Algorithm B. We show that the vertices can be thrown out in the order $\pi(1), 
\pi(2), \ldots, \pi(n)$.  As $\pi(1) = 1$, the vertex $\pi(1)$ clearly can be thrown out.  Suppose $\pi(1), \ldots, \pi(k-1)$ can be thrown out, and 
consider $\pi(k)$.

If $f(\pi(k)) = \infty$, by Case (1) of step 1, $\pi(k)$ is the smallest vertex not in $\{\pi(1), \dots, \pi(k-1)\}$. Thus it can be thrown out.

If $f(\pi(k)) \neq \infty$, there is an edge of the forest $F$ connecting $\pi(k)$ to a vertex $w$ in $\{\pi(1), \dots, \pi(k-1)\}$. Suppose $w = 
\pi(t)$ where $t < k$.  By definition of $f$, there are exactly $f(\pi(k))$ edges connecting $\pi(k)$ to the set $\{ \pi(1), \dots, \pi(t-1)\}$.  
Hence $f(\pi(k)) < outdeg_{\{\pi(k), \dots, \pi(n)\}}(\pi(k))$. Thus $\pi(k)$ can be thrown out as well.

By induction the vertices of $G$ can be thrown out in the order $\pi(1), \pi(2), \ldots, \pi(n)$.

\end{proof}
 
Define $\Psi_{\gamma, G}: \fG \rightarrow \mpf$ by letting $\Psi_{\gamma, G}(f)=f_F$.  Now we show that $\Phi=\Phi_{\gamma,G}$ and 
$\Psi=\Psi_{\gamma,G}$ are inverses of each other.

\begin{thm} \label{bijection} $\Psi(\Phi(f))=f$ for any $f \in \mpf$ and $\Phi(\Psi(F))=F$ for any $F \in \fG$.

\end{thm}

\begin{proof}

First, if $f \in \mpf$ and $F=\Phi(f)$, then by Prop. \ref{root} the roots of $f$ are exactly the minimal vertices in each tree component of $F$. 
Those in turn are roots for $\Psi(F)$.  In applying algorithm B to $F$, we note that the order $\pi=v_1v_2\dots v_n$ is exactly the order in which 
vertices of $G$ will be processed when running algorithm A on $f$. That is, $P_i=\{v_1, \dots, v_i\}$, and $v_{i+1}$ is not a root of $f$, then 
$Q_i$ is the set of vertices which are adjacent (via edges in $F$) to those in $P_i$. By the construction of algorithm A, a vertex $w$ is found by 
$v$ if and only if there are $f(w)$ many edges connecting $w$ to vertices that are processed before $v$, or equivalently, to vertices $u$ with 
$\pi^{-1}(u) < \pi^{-1}(v)$. Since in $\Phi(f)$, $v=w^p$, we have $\Psi(\Phi(f))=f$.

Conversely, we prove that $\Phi(\Psi(F))=F$ by showing that $\Phi(\Psi(F))$ and $F$ have the same set of edges.  First note that the minimal 
vertices of the tree components of $F$ are exactly the roots of $f=\Psi(F)$, which then are the minimal elements of trees in $\Phi(f)$.  Edges of 
$F$ are of the form $\{v, v^p\}$, where $v$ is not a minimal vertex in its tree component.  We now show that when applying algorithm A to $\Psi(F)$, 
vertex $v$ is found by $v^p$. Note that $f(v)= |\{ v_j | ( v,v_j ) \in E(G), \pi^{-1}(v_j) < \pi^{-1}(v^p)\}|$. In the implementation of algorithm 
A, the {\em valuation} on $v$ drops by 1 for each adjacent vertex that is processed before $v$.  When it is $v^p$'s turn to be processed, $val_i(v)$ 
drops from $0$ to $-1$. Thus $v^p$ finds $v$, and $\{v, v^p\}$ is an edge of $\Phi(\Psi(F))$. 

\end{proof}

Since the roots of the \gmp \ correspond exactly to the minimal vertices in the tree components of the corresponding forest, in the following we 
will refer to those vertices as \emph{roots of the forest}.

\section{Examples of the bijections}

The bijections $\Phi_{\gamma, G}$ and $\Psi_{\gamma, G}$, as defined above via algorithms A and B, allow a good deal of freedom in implementation.  
In algorithm A, as long as $\gamma$ is well-defined at every iteration of Step 2, one can obtain $val_{i+1}$, $P_{i+1}$, $Q_{i+1}$ and $F_{i+1}$ and 
proceed. Recall that $\gamma$ is a function from $\prod$, the set of ordered pairs $(F, W)$, to $V(G)$ such that $\gamma(F, W) \in W$, where $F$ is 
any sub-forest of $G$ (not necessarily spanning) and $W$ is a non-empty subset of $Leaf(F)$ or consists of an isolated point of $F$.

When restricting to $G$-parking functions, (i.e., $G$-multiparking functions with only one root), the descriptions of the bijections $\Phi$ and 
$\Psi$ are basically the same as the ones given by Chebikin and Pylyavskyy \cite{ChPy:2005}, where the corresponding sub-structures in $G$ are 
spanning trees. However our family of bijections, each defined on a choice function $\gamma$, is more general than the ones in \cite{ChPy:2005}, 
which rely on a {\em proper set of tree orders}. A proper set of tree orders is a set $\Pi(G)=\{\pi(T): T \text{ is a subtree of } G \}$ of linear 
orders on the vertices of $T$, such that for any $v \in T$, $v <_{\pi(T)} v^p$, and if $T'$ is a subtree of $T$ containing the least vertex, 
$\pi(T')$ is a suborder of $\pi(T)$.  Our algorithms do not require there to be a linear order on the vertices of each subtree. In fact, for a 
spanning tree $T$ of a connected graph $G$, the proper tree order $\pi(T)$, if it exists, must be the same as the one defined in Step 1 of algorithm 
B. But in general, for two spanning trees $T$ and $T'$ with a common subtree $t$, the restrictions of $\pi(T)$ and $\pi(T')$ to vertices of $t$ may 
not agree. Hence in general the choice function cannot be described in terms of proper sets of tree orders. In addition, our description of the map 
$\Phi$, in terms of a dynamic process, provides a much clearer way to understand the bijection, and leads to a natural classification of the edges 
of $G$ which plays an important role in connection with the Tutte polynomial (c.f.\S4).

\begin{figure}[ht] 
\[
\begin{picture}(0,0)
\put(-95,-137){\includegraphics[width=7cm]{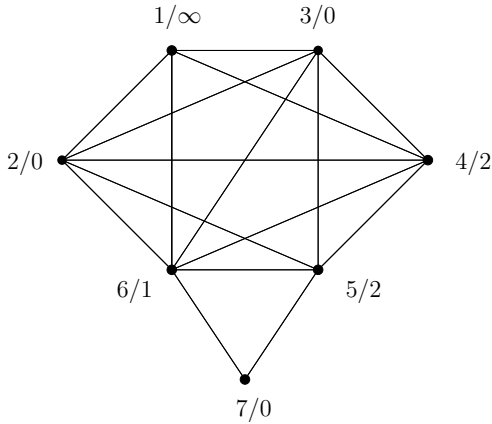}}
\end{picture}
\]
\mbox{}\vspace{1.5in} 
\caption{A graph and a multiparking function.}\label{graph-example}
\end{figure}

Different choice functions $\gamma$ will induce different bijections between $\mpf$ and $\fG$. In this section we give several examples of choice 
functions that have combinatorial significance. In Example 1 we explain how to translate a proper set of tree orders into a choice function.  Hence 
the family of bijections defined in \cite{ChPy:2005} can be viewed as a subfamily of our bijections restricted to $G$-parking functions.  The next 
three examples have appeared in \cite{ChPy:2005}. We list them here for their combinatorial significance.  Example 5 is the combination of 
breadth-first search with the $Q$-sets equipped with certain data structures.  It is the one used to establish connections with Tutte polynomial in 
\S4. The last example illustrates a case where $\gamma$ cannot be expressed as a proper set of tree orders.  We illustrate the corresponding map 
$\Phi_{\gamma, G}$ for examples 2--6 on the graph $G$ in Figure \ref{graph-example}.  A $G$-multiparking function $f$ is indicated by ``$i/f(i)$" on 
vertices, where $i$ is the vertex label.

In each example, we will show the resulting spanning forest by darkened edges in $G$.  Again each vertex will be labeled by a pair $i/j$, where $i$ 
is the vertex labels, and $j=val_n(i)$, where $n=7$.  Beneath that, a table will record the sets $Q_t$ and $P_t$ for each time $t$. In each $Q_t$, 
the vertex listed first is the next to be processed.

\noindent \textbf{Example 1.} $\gamma$ with a proper set of tree orders. \\ We define the choice function that corresponds to a proper set of tree 
orders. Here we should generalize to the proper set of forest orders, i.e., a set of orders $\pi(F)$, defined on the set of vertices for each 
subforest $F$ of $G$, such that for any $v \in F$, $v <_{\pi(F)} v^p$, and if $F'$ is a subforest of $F$ with the same minimal vertex in each tree 
component, $\pi(F')$ is a suborder of $\pi(F)$. In this case, define $\gamma(F,W)=v$ where $v$ is the minimal element in $W$ under the order 
$\pi(F)$. Examples 2--4 are special cases of this kind.

\noindent {\bf Example 2.} $\gamma$ with a given vertex ranking. \\ Given a vertex ranking $\sigma \in S_n$ define $\gamma_\sigma(F, W):=v$, where 
$v$ is the vertex in $W$ with minimal ranking. In particular, if $\sigma$ is the identity permutation, then the vertex processing order is the {\em 
vertex-adding order} of \cite{ChPy:2005}.  In this case, in Step 2 of algorithm A, we choose $v$ to be the least vertex in $Q_{i-1}$ and process it 
at time $i$. The output of algorithm A is

\begin{figure}[ht]
\[
\begin{picture}(0,0)
\put(-95,-130){\includegraphics[width=7cm]{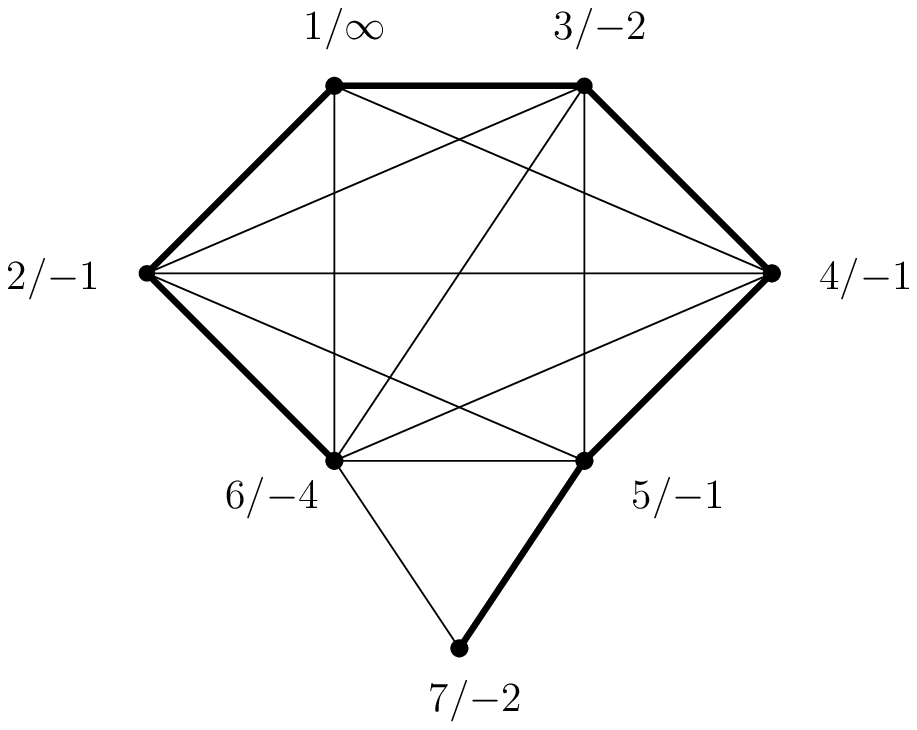}}
\end{picture}
\]
\mbox{} \vskip 1.4in
\end{figure}

\noindent The $Q_i$ and $P_i$ for this instance are as follows.
\begin{center}
\begin{tabular}{|c|c|c|c|c|c|c|c|c|}
\hline
\hspace{.2cm} t\hspace{.2cm} &
\hspace{.2cm} 0\hspace{.2cm} &
\hspace{.2cm} 1\hspace{.2cm} &
\hspace{.2cm} 2\hspace{.2cm} &
\hspace{.2cm} 3\hspace{.2cm} &
\hspace{.2cm} 4\hspace{.2cm} &
\hspace{.2cm} 5\hspace{.2cm} &
\hspace{.2cm} 6\hspace{.2cm} &
\hspace{.2cm} 7\hspace{.2cm} \\ \hline
$Q_t$ & \{1\} & \{2,3\} & \{3,6\} & \{4,6\} & \{5,6\} & \{6,7\} & \{7\}& $\emptyset$ \\ \hline
$P_t$ & $\emptyset$ & \{1\} & \{1,2\} & \{1,2,3\} & \{1,2,3,4\} & \{1,2,3,4,5\} & \{1,2,3,4,5,6\}  & \{1,2,3,4,5,6,7\} \\ \hline
\end{tabular}
\end{center}

\noindent{\bf Example 3.} $\gamma$ with \emph{depth-first search order}.\\

The depth-first search order is the order in which vertices of a forest are visited when performing the depth-first search, which is also known as 
the preorder traversal. Given a forest $F$ with tree components $T_1, T_2, \dots, T_k$, where $1=r_1< r_2< \cdots <r_k$ are the corresponding roots, 
the order $<_{df}$ is defined as follows. (1) For any $v \in T_i$, $w \in T_j$ and $i<j$, $v <_{df} w$. (2) For any $v \neq r_i$, $v^p <_{df} v$. 
(3) If $v^p=w^p$ and $v<w$, $v <_{df} w$. (4) For any $v$, let $F[v]$ be the subtree of $F$ rooted at $v$. If $v \in F[v']$, $w \in F[w']$ and $v' 
<_{df} w'$, then $v <_{df} w$. For example, the depth-first search order on the below tree is $1 <_{df} 2 <_{df} 3 <_{df} 6 <_{df} 4 <_{df} 5$.

\begin{figure}[ht]
\[
\begin{picture}(0,0)
\put(-55,-45){\includegraphics[width=4cm]{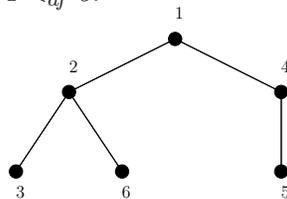}}
\end{picture}
\]
\vskip 0.3 in
\caption{A tree with 6 vertices.} \label{tree-example} 
\end{figure}

The choice function $\gamma_{df}$ with depth-first search order is then defined as $\gamma_{df}(F, W)=v$ where $v$ is the minimal element of $W$ 
under the depth-first search order $<_{df}$ of $F$. Here is the output of algorithm A with the choice function $\gamma_{df}$ on the example in 
Figure \ref{graph-example}.

\begin{figure}[ht]
\[
\begin{picture}(0,0)
\put(-95,-137){\includegraphics[width=7cm]{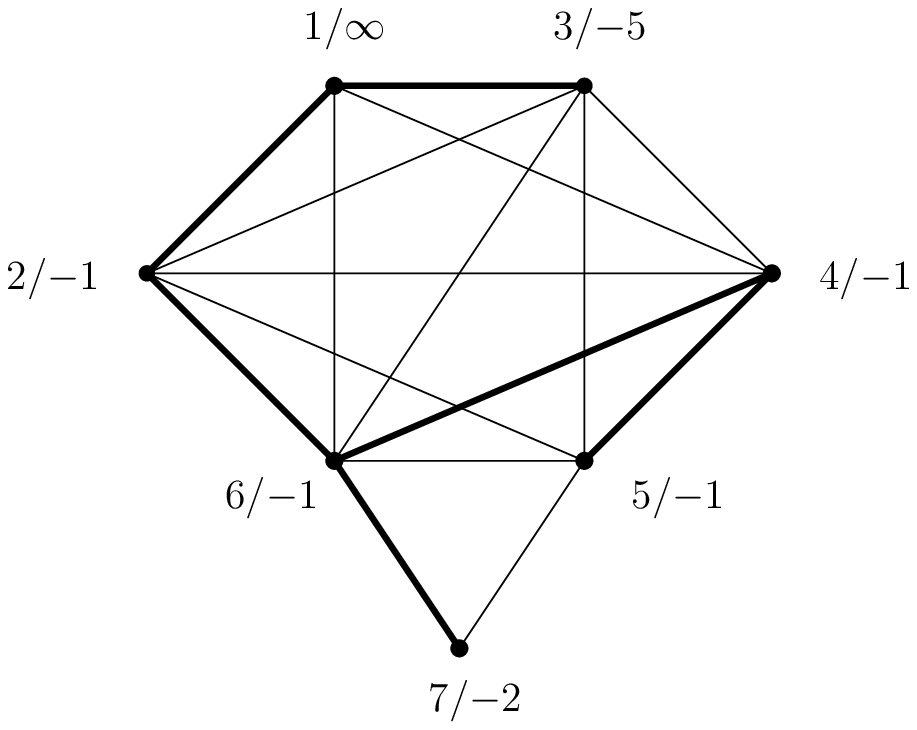}}
\end{picture} 
\]
\mbox{} \vskip 1.4in
\end{figure}

\noindent The $Q_i$ and $P_i$ for this instance are as follows.
\begin{center}
\begin{tabular}{|c|c|c|c|c|c|c|c|c|}
\hline
\hspace{.2cm} t\hspace{.2cm} &
\hspace{.2cm} 0\hspace{.2cm} &
\hspace{.2cm} 1\hspace{.2cm} &
\hspace{.2cm} 2\hspace{.2cm} &
\hspace{.2cm} 3\hspace{.2cm} &
\hspace{.2cm} 4\hspace{.2cm} &
\hspace{.2cm} 5\hspace{.2cm} &
\hspace{.2cm} 6\hspace{.2cm} &
\hspace{.2cm} 7\hspace{.2cm} \\ \hline
$Q_t$ & \{1\} & \{2,3\} & \{6,3\} & \{4,3,7\} & \{5,3,7\} & \{7,3\} & \{3\} & $\emptyset$ \\ \hline
$P_t$ & $\emptyset$ & \{1\} & \{1,2\} & \{1,2,6\} & \{1,2,4,6\} & \{1,2,4,5,6\} & \{1,2,4,5,6,7\} & \{1,2,3,4,5,6,7\} \\ \hline
\end{tabular}
\end{center}

\noindent {\bf Example 4.} $\gamma$ with {\em breadth-first search  order}. \\ 
Breadth-first search is another commonly used tree traversal in computer science. Given a forest $F$, whose tree components are $T_i$ with roots 
$r_i$, ($1\leq i \leq k$), and $1=r_1< r_2< \cdots <r_k$, the order $<_{bf}$ is defined as follows. (1) For any $v \in T_i$, $w \in T_j$ and $i<j$, 
$v <_{{bf}} w$. (2) Within tree $T_i$, for each $v \in T_i$, let \emph{height} $h_{T_i}(v)$ of $v$ be the number of edges in the unique path from 
$v$ to the root $r_i$. We set $v <_{{bf}} w$ if $h_{T_i}(v) < h_{T_i}(w)$, or else if $h_{T_i}(v) = h_{T_i}(w)$ and $v < w$.  For example, the the 
breadth-first search order for the tree in Figure \ref{tree-example} is $1 <_{bf} 2 <_{bf} 4 <_{bf} 3 <_{bf} 5 <_{bf} 6$.

The choice function $\gamma_{bf}$ with breadth-first search order is defined as $\gamma_{bf}(F, W)=v$ where $v$ is the minimal element of $W$ under 
the breadth-first search order $<_{bf}$ of $F$. Here is the output of algorithm A with the choice function $\gamma_{bf}$ on the example in Figure 
\ref{graph-example}.

\begin{figure}[ht]
\[
\begin{picture}(0,0)
\put(-95,-137){\includegraphics[width=7cm]{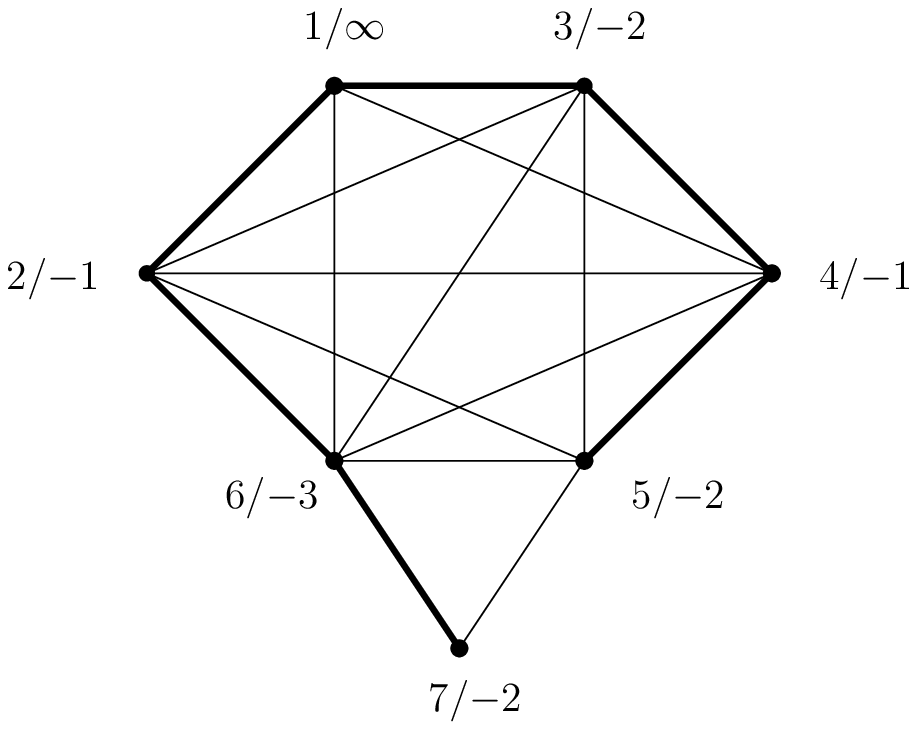}}
\end{picture}
\]
\mbox{} 
\end{figure}

\pagebreak
\noindent The $Q_i$ and $P_i$ for this instance are as follows.
\begin{center}
\begin{tabular}{|c|c|c|c|c|c|c|c|c|}
\hline
\hspace{.2cm} t\hspace{.2cm} &
\hspace{.2cm} 0\hspace{.2cm} &
\hspace{.2cm} 1\hspace{.2cm} &
\hspace{.2cm} 2\hspace{.2cm} &
\hspace{.2cm} 3\hspace{.2cm} &
\hspace{.2cm} 4\hspace{.2cm} &
\hspace{.2cm} 5\hspace{.2cm} &
\hspace{.2cm} 6\hspace{.2cm} &
\hspace{.2cm} 7\hspace{.2cm} \\ \hline
$Q_t$ & \{1\}& \{2,3\} & \{3,6\} & \{4,6\} & \{6,5\} & \{5,7\} & \{7\}& $\emptyset$ \\ \hline
$P_t$ & $\emptyset$ & \{1\} & \{1,2\} & \{1,2,3\} & \{1,2,3,4\}  &
  \{1,2,3,4,6\}  & \{1,2,3,4,5,6\} & \{1,2,3,4,5,6,7\} \\ \hline
\end{tabular}
\end{center}

\noindent{\bf Example 5.} Breadth-first search with a data structure on $Q_i$. 

In this case, new vertices enter the set $Q_i$ in a certain order, and some intrinsic data structure on $Q_i$ decides which vertex of $Q_i$ is to be 
processed in the next step. A typical example is that of breadth-first search with a queue, in which case each $Q_i$ is an ordered set, (i.e., the 
stage of a queue at time $i$). New vertices enter $Q_i$ in numerical order, and $\gamma$ chooses the vertex that entered the queue earliest.

This example can also be defined by a modified breadth-first search order, which we call \emph{breadth-first order with a queue}, and denote by 
$<_{bf, q}$. Given a forest $F$, whose tree components are $T_i$ with root $r_i$, ($1\leq i \leq k$), and $1=r_1< r_2< \cdots <r_k$, the order 
$<_{bf,q}$ is defined as follows. (1) For any $v \in T_i$, $w \in T_j$ and $i<j$, $v <_{bf,q} w$. (2) Within tree $T_i$, the root $r_i$ is minimal 
under $<_{bf,q}$. (3) $v <_{{bf,q}} w$ if $v^p <_{{bf,q}} w^p$.  (4) If $v^p=w^p$ and $v < w$, $v <_{{bf,q}} w$. For example, the breadth-first 
search order with a queue for the tree in Figure \ref{tree-example} is $1 <_{bf,q} 2 <_{bf,q} 4 <_{bf,q} 3 <_{bf,q} 6 <_{bf,q} 5$.

The choice function $\gamma$ associated with this order is denoted by $\gamma_{bf,q}$, and is used in \S4. The following is the output of algorithm 
A with $\gamma_{bf,q}$ on the graph in Figure \ref{graph-example}.

\begin{figure}[ht]
\[
\begin{picture}(0,0)
\put(-95,-137){\includegraphics[width=7cm]{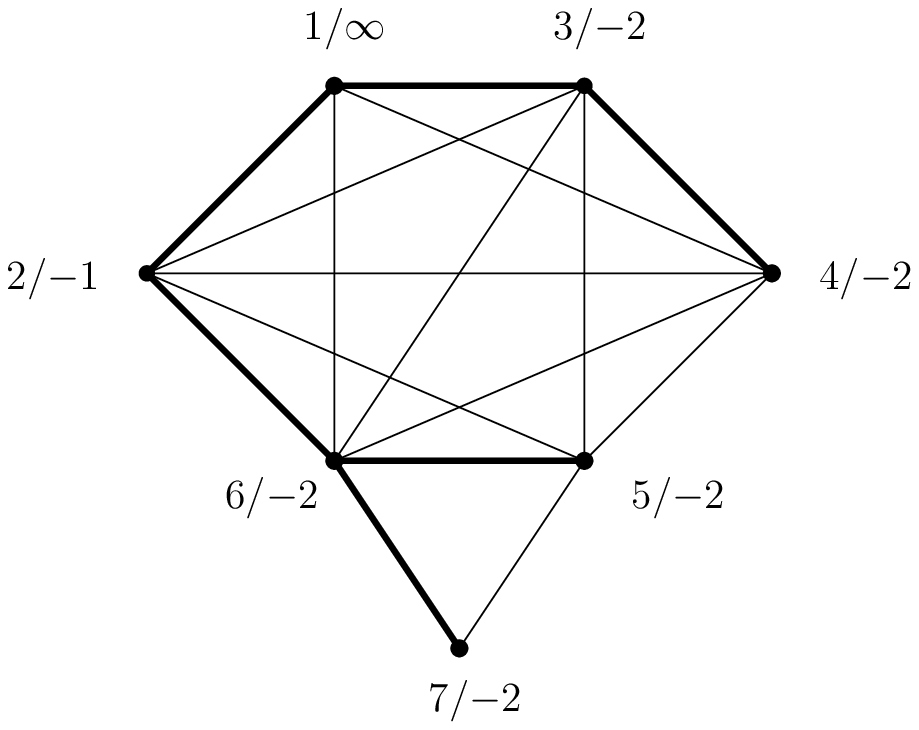}}
\end{picture}
\]
\mbox{} \vskip 1.5 in
\end{figure}

\noindent The $Q_i$ and $P_i$ for this instance are as follows, where each $Q_i$ is an ordered set, and the first element in $Q_i$ is the next one 
to be processed.
   
\begin{center}
\begin{tabular}{|c|c|c|c|c|c|c|c|c|}
\hline
\hspace{.2cm} t\hspace{.2cm} &
\hspace{.2cm} 0\hspace{.2cm} &
\hspace{.2cm} 1\hspace{.2cm} &
\hspace{.2cm} 2\hspace{.2cm} &
\hspace{.2cm} 3\hspace{.2cm} &
\hspace{.2cm} 4\hspace{.2cm} &
\hspace{.2cm} 5\hspace{.2cm} &
\hspace{.2cm} 6\hspace{.2cm} &
\hspace{.2cm} 7\hspace{.2cm} \\ \hline
$Q_t$ & (1) & (2,3) & (3,6) & (6,4) & (4,5,7) & (5,7) & (7) & $\emptyset$ \\ \hline
$P_t$ & $\emptyset$ & \{1\} & \{1,2\} & \{1,2,3\} & \{1,2,3,6\} & \{1,2,3,4,6\} & \{1,2,3,4,5,6\}  & \{1,2,3,4,5,6,7\} \\ \hline
\end{tabular}
\end{center}

Another typical structure is to let $Q_i$ be the stage of a stack at time $i$, that is, it pops out the vertex that last entered. We can also 
combine the other vertex orders with a queue or stack for the $Q$-sets.

\noindent {\bf Example 6.} A choice function $\gamma$ that cannot be defined by a proper set of tree orders.

Let 
\begin{eqnarray*} 
\gamma(F, W)=\left\{\begin{array}{ll} 
  x &  \text{ if } W=\{x\}, \\
  \text{the second minimal vertex of $W$}, & \text{ if } |W| \geq 2.
\end{array} \right. 
\end{eqnarray*} 

Then the order on the left tree is $156342$, and the one on the right tree is $153462$, which do not agree on the subtree consisting of vertices 
$1356$. Hence it can not be defined via a proper set of tree orders.

\begin{figure}[ht] 
\centerline{\epsfxsize=2.5in \epsfbox{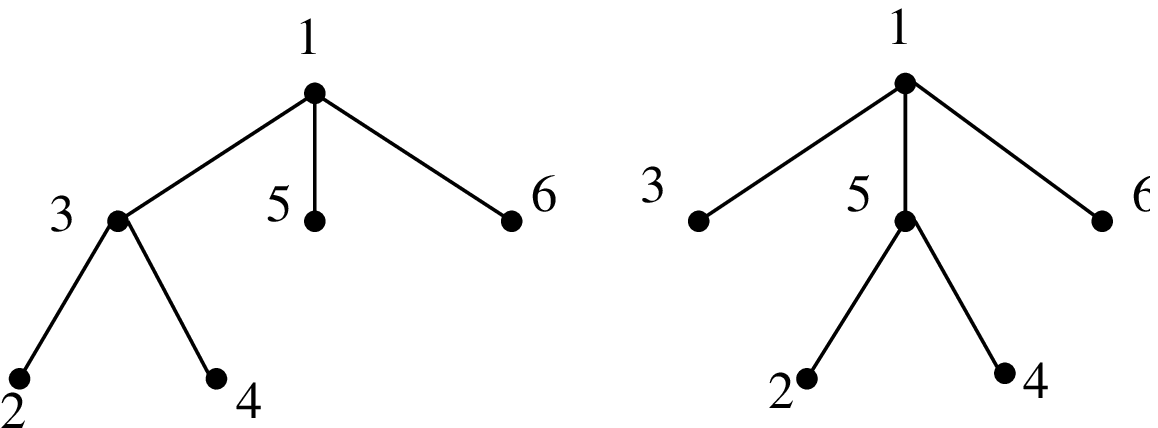}}
\end{figure} 

\textsc{Remark}.  Note that the bijection given in \cite{ChPy:2005} applies to a general directed graph. We explain how our algorithms could be 
slightly modified to apply to that case, too. Let $D$ denote a general directed graph on $[n]$. An oriented spanning forest $F$ of $D$ is a subgraph 
of $G$ such that (1) the edges, when ignoring the orientation, do not form a cycle, and (2) for each $v$ in a tree component with minimal vertex 
$r$, there is a unique directed path from $v$ to $r$. Again we denote by $v^p$ the vertex lying on the directed path from $v$ to $r$ with $(v, 
v^p)\in E(D)$. We say that the minimal vertex in each tree component of $F$ is a root of $F$.

The definition of a $D$-multiparking function is the same as that of the undirected multiparking function, except that $outdeg_U(i)$ is the number 
of edges going from $i$ to vertices not in $U$.  Again we say the vertices $v$ with $g(v) = \infty$ are the roots of the $D$-multiparking function 
$g$.

Spanning forests do not contain loops, as loops are a trivial kind of cycle. Hence we can assume $D$ is loopless without loss of generality. We 
allow $D$ to have multiple edges. But to distinguish between multiple edges of $D$, we fix a total order on the set of edges going from $i$ to $j$, 
for each $i\neq j$.

The maps $\Phi$ and $\Psi$ can be modified accordingly to give a bijection between the set of $D$-multiparking functions to the set of oriented 
spanning forests of $D$, which carry the roots of multiparking functions to the roots of spanning forests.  The only modifications we need to make 
are:

\noindent \textbf{For Algorithm A.}

In Step 3, lower the value of $w$ by 1 for each directed edge from $w$ to $v$ if $w$ is not a root. Load $w$ to $Q_i$ whenever $val_i(w)<0$.  For 
each such $w$, add the $(k+1)$-st edge between from $w$ to $v$ if $val_{i-1}(w)=k \geq 0$.

\noindent \textbf{For Algorithm B.}

In Step 1. Let $W$ be the set $\{v \notin V_i: \exists w \in V_i$ such that $(v,w)\in E(D)\}$.
 
In Step 2. For any nonroot vertex $v$ lying in the tree with root $r_v$, and $v \neq r_v$, set $f(v)$ to be $k+ |\{ v_j | ( v, v_j ) \in E(G), 
\pi^{-1}(v_i) < \pi^{-1}(v^p)\}|$ if the edge $(v,v^p)$ in $F$ is the $(k+1)$-st edge in the set of edges from $v$ to $v^p$ in $D$.

\vspace{.2cm}

These modified algorithms for directed graphs cover the case for an undirected graph $G$, provided that one views $G$ as a digraph, where each edge 
$\{u,v\}$ of $G$ is replaced with two directed edges $(u,v)$ and $(v,u)$.

The notion of $G$-parking functions, as proposed in \cite{PostnikovShapiro:2004}, is closely related to the critical configurations of the 
chip-firing games, (also known as sand-pile models).  The generalization of chip-firing games with multiple sources is given in \cite{Ellis:2002}, 
where they are called \emph{Dirichlet games}. In \cite{Kostic:2006} the first author shows that modified $G$-multiparking functions (which in 
Condition {\bf (A)}, instead of requiring that $i=min(U)$, one requires $i$ to belong to a prefixed subset of vertices), are the corresponding 
counterpart for critical configurations of the Dirichlet games. In fact, both are in one-to-one correspondence with the set of rooted spanning 
forests of $G$, as well as a set of objects called descending multitraversals on $G$.

%
%
\section{External activity and the Tutte polynomial} 

\subsection{$F$-redundant edges}
A forest $F$ on $[n]$ may appear as a subgraph of different graphs, and a vertex function $f$ may be a \gmp \ for different graphs. In this section 
we characterize the set of graphs which share the same pair $(F, f)$. Again let $G$ be a simple graph on $[n]$, and fix a choice function $\gamma$. 
For a spanning forest $F$ of $G$, let $f=\Psi_{\gamma,G}(F)$.  We say an edge $e$ of $G-F$ is \emph{$F$-redundant} if $\Psi_{\gamma, G -\{ e \}}(F) 
= f$. Note that we only need to use the value of $\gamma$ on $(F', W)$ where $F'$ is a sub-forest of $F$. Hence $\Psi_{\gamma, G -\{ e \}}(F)$ is 
well-defined.

Let $\pi$ be the order defined in Step 1 of Algorithm B. Note that $\pi$ only depends on $F$, not the underlying graph $G$. Recall that $v^p$ 
denotes the parent vertex of vertex $v$ in some spanning forest.  We have the following proposition.

\begin{prop}~\label{edgecharacterization} 

An edge $e = \{v,w\}$ of $G$ is $F$-redundant if and only if
$e$ is one of the following types:

\begin{enumerate}

\item Both $v$ and $w$ are roots of $F$.

\item $v$ is a root and $w$ is a non-root of $F$, and $\pi^{-1}(w) < \pi^{-1}(v)$.

\item $v$ and $w$ are non-roots and $\pi^{-1}(v^p)< \pi^{-1}(w) < \pi^{-1}(v)$. In this case $v$ and $w$ must lie in the same tree of $F$.

\end{enumerate} 

\end{prop}

\begin{proof}

We first show that each edge of the above three types are $F$-redundant. Since for any root $r$ of the forest $F$, $f(r) = \infty$, the edges of the 
first two types play no role in defining the function $f$.  And clearly those edges are not in $F$. Hence they are $F$-redundant.

For edge $(v,w)$ of type 3, clearly it cannot be an edge of $F$.  Since $f(v)= \# \{ v_j | ( v, v_j ) \in E(G), \pi^{-1}(v_j) < \pi^{-1}(v^p)\}$, 
and $\pi^{-1}(w) > \pi^{-1}(v^p)$, removing the edge $\{v,w\}$ would not change the value of $f(v)$.  This edges has no contribution in defining 
$f(u)$ for any other vertex $u$. Hence it is $F$-redundant.

For the converse, suppose that $e=\{v,w\}$ is not one of the three type. Assume $w$ is processed before $v$ in $\pi$. Then $v$ is not a root, and 
$w$ appears before $v^p$. Then removing the edge $e$ will change the value of $f(v)$. Hence it is not $F$-redundant. 

\end{proof}

Let $R_1(G;F)$, $R_2(G;F)$, and $R_3(G;F)$ denote the sets of $F$-redundant edges of types $1$, $2$, and $3$, respectively.  Among them, $R_3(G; F)$ 
is the most interesting one, as $R_1(G;F)$ and $R_2(G;F)$ are a consequence of the requirement that $f(r)=\infty$ for any root $r$. Let $R(G;F)$ be 
the union of these three sets. Clearly the $F$-redundant edges are mutually independent, and can be removed one by one without changing the 
corresponding \gmp. Hence 

\begin{thm}~\label{preimage} Let $H$ be a subgraph of $G$ with $V(H)=V(G)$.  Then $\Psi_{\gamma, G}(F)=\Psi_{\gamma, 
H}(F)$ if and only if $G - R(G;F) \subseteq H \subseteq G$. 

\end{thm}

\subsection{A classification of the edges of $G$} 

The notion of $F$-redundancy allows us to classify the edges of a graph in terms of the algorithm A.  Roughly speaking, the edges of any graph can 
be thought of as either lowering $val(v)$ for some $v$ to $0$, being in the forest, or being $F$-redundant.  Explicitly, we have

\begin{prop}~\label{edgecategorization} 

Let $f$ be a $G$-multiparking function and let $F = \Phi(f)$.  Then
\[ 
|E(G)| = \bigg( \sum_{v : f(v) \neq \infty} f(v) \bigg) + |E(F)| +
 |R(G;F)|. 
\]

\end{prop}

\begin{proof}

For each non-root vertex $v$, the number of different values that $val_i(v)$ takes on during the execution of algorithm $A$ is $f(v)+1+n_v$, where 
$n_v=-val_n(v)$. At the beginning, $val_0(v) = f(v)$. The value $val_i(v)$ then is lowered by one whenever there is a vertex $w$ which is adjacent 
to $v$ and processed before $v^p$. When $v^p$ is being processed, $val_i(v)=-1$, and the edge $\{v^p, v\}$ contributes to the forest $F$. Afterward, 
the value of $val_i(v)$ decreases by 1 for each $F$-redundant edge $\{u,v\}$ with $\pi^{-1}(u) < \pi^{-1}(v)$.  Summing over all non-root vertices 
gives \begin{eqnarray*} \sum_{v : f(v) \neq \infty} deg_{<_{\pi}}(v) &=& \sum_{v : f(v) \neq \infty} f(v) + |E(F)| + \sum_{v:f(v)\neq \infty} n_v, 
\end{eqnarray*} where $deg_{<_{\pi}}(v)=|\{ \{w, v\} \in E(G)| \pi^{-1}(w) < \pi^{-1}(v)\}|$.

The edges that lower $val(v)$ below $-1$ are exactly the $F$-redundant edges of type (3) in Prop. \ref{edgecharacterization}, hence 
$\sum_{v:f(v)\neq \infty} n_v =|R_3(G;F)|$. On the other hand, $\sum_{v : f(v) \neq \infty} deg_{<_{\pi}}(v)$ is exactly $|E(G)|
 -|R_1(G;F)|-|R_2(G;F)|$. The claim follows from the fact that the sets $R_1(G;F), R_2(G;F)$, and $R_3(G;F)$ are mutually exclusive.

\end{proof} 

One notes that for roots of $f$ and $F=\Phi(f)$, $|R_1(G;F)|+|R_2(G;F)|$ is exactly $\sum_{root\ v} deg_{<_{\pi}}(v)$, where $\pi$ is the processing 
order in algorithm A. But it is not necessary to run the full algorithm A to compute $|R_1(G;F)|+|R_2(G;F)|$. Instead, we can apply the burning 
algorithm in a greedy way to find an ordering $\pi'=v_1'v_2' \cdots v_n'$ on $V(G)$: Let $v_1'=1$. After determining $v_1', \dots, v_{i-1}'$, if 
$V_i=V(G)-\{v_1, \dots, v_{i-1}'\}$ has a well-behaved vertex, let $v_i'$ be one of them; otherwise, let $v_i'$ be the minimal vertex of $V_i$, 
(which has to be a root.)

$\pi'$ may not be the same as $\pi$, but they have the following properties:

\begin{enumerate}

\item Let $r_1< r_2< \cdots < r_k$ be the roots of $f$. Then $r_1, r_2, \dots, r_k$ appear in the same positions in both $\pi$ and $\pi'$.

\item The set of vertices lying between $r_i$ and $r_{i+1}$ are the same in $\pi$ and $\pi'$. In fact, they are the vertices of the tree $T_i$ with 
root $r_i$ in $F=\Phi(f)$.

\end{enumerate} 

It follows that for any root vertex $v$, $deg_{<_\pi}(v)=deg_{<_{\pi'}}(v)$.  The value of $deg_{<_\pi}(v)$ ($v$ root) can be characterized by a 
global description: Let $\mathcal{U}_v$ be the collection of subsets $U$ of $V(G)$ such that $v=\min(U)$, and $U$ does not have a well-behaved 
vertex. $\mathcal{U}_v$ is nonempty for a root $v$ since $U=\{v\}$ is such a set.  Then 

$$ 
deg_{<_\pi}(v)=\min_{U \in \: \mathcal{U}_v}outdeg_U(v). 
$$

We call $deg_{<_\pi}(v)$ the \emph{record} of the root $v$, and denote it by $rec(v)$. Then 

$$
|R_1(G;F)|+|R_2(G;F)| = \sum_{root\ v} deg_{<_{\pi'}}(v) =\sum_{root \ v} rec(v)
$$

is the \emph{total root records}. Let $Rec(f)=|R_1(G;F)|+|R_2(G;F)|$.  It is the number of $F$-redundant edges adjacent to a root. By the above 
greedy burning algorithm, the total root records $Rec(f)$ can be computed in linear time.

\subsection{A new expression for Tutte polynomial} 

In this subsection we relate $G$-multiparking functions to the Tutte polynomial $t_G(x,y)$ of $G$.  We follow the presentation of 
\cite{GesselSagan:1996} for the definition of Tutte polynomial and its basic properties. Although the theory works for general graphs with 
multiedges, we assume $G$ is a simple connected graph to simplify the discussion. There is no loss of generality by assuming connectedness, since 
for a disconnected graph, $t_G(x,y)$ is just the product of the Tutte polynomials of the components of $G$. We restrict ourselves to connected 
graphs to avoid any possible confusion when we consider their spanning forests.  The modification when $G$ has multiple edges is explained at the 
end of \S3.

Suppose we are given $G$ and a total ordering of its edges. Consider a spanning tree $T$ of $G$. An edge $e \in G-T$ is \emph{externally active} if 
it is the largest edge in the unique cycle contained in $T \cup e$.  We let $$ \mathcal{EA}(T)=\text{set of externally active edges in } T $$ and 
$ea(T)=|\mathcal{EA}(T)|$. An edge $e \in T$ is {\em internally active} if it is the largest edge in the unique cocycle contained in $(G-T)\cup e$. 
We let $$ \mathcal{IA}(T)=\text{set of internally active edges in } T $$ and $ia(T)=|\mathcal{IA}(T)|$. Tutte \cite{Tutte:1953} then defined his 
polynomial as \begin{eqnarray} t_G(x,y) =\sum_{T \subset G} x^{ia(T)} y^{ea(T)}, \end{eqnarray} where the sum is over all spanning trees $T$ of $G$. 
Tutte showed that $t_G$ is well-defined, i.e., independent of the total ordering of the edges of $G$. Henceforth, we will not assume that the edges 
of $G$ are ordered.

Let $H$ be a (spanning) subgraph of $G$. Denote by $c(H)$ the number of components of $H$. Define two invariants associated with $H$ as 

\begin{eqnarray} 
\sigma(H)=c(H)-1, \qquad \sigma^*(H)=|E(H)|-|V(G)|+c(H). 
\end{eqnarray}

The following identity is well-known, for example, see \cite{Biggs:1993}.
\begin{thm} 
\begin{eqnarray} \label{Tutte-2} 
t_G(1+x,1+y) = \sum_{H \subseteq G} x^{\sigma(H)} y^{\sigma^*(H)},
\end{eqnarray} 
where the sum is over all spanning subgraphs $H$ of $G$.  
\end{thm} 

Recall that the \emph{breadth-first search} (BFS) is an algorithm that gives a spanning forest in the graph $H$. Assume $V(G)=[n]$.  We will use our 
favorite description to express the BFS as a queue $Q$ that starts at the least vertex $1$. This description was first introduced in 
\cite{Spencer:1997} to develop an exact formula for the number of labeled connected graphs on $[n]$ with a fixed number of edges, and was used by 
the second author in \cite{Yan:2001} to reveal the connection between the classical parking functions (resp. $k$-parking functions) and the complete 
graph (resp. multicolored graphs).

\begin{figure}[ht]
\centerline{\epsfxsize=5in \epsfbox{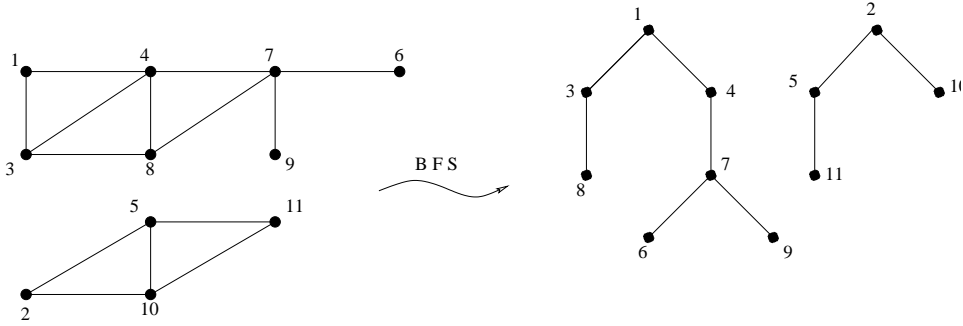}}
\caption{Spanning forest found by BFS.} \label{BFS_on_H}
\end{figure}

Given a subgraph $H$ of $G$ with $V(H)=V(G)=[n]$, we construct a queue $Q$. At time $0$, $Q$ contains only the vertex $1$. At each stage we take the 
vertex $x$ at the head of the queue, remove $x$ from the queue, and add all unvisited neighbors $u_1, \dots, u_{t_x}$ of $x$ to the queue, in 
numerical order.  We will call this operation ``processing $x$''.  If the queue becomes empty, add the least unvisited vertex to $Q$.  The output 
$F$ is the forest whose edge set consists of all edges of the form $\{x, u_i\}$ for $i=1, \dots t_x$.  We will denote this output as $F=BFS(H)$.  
Figure \ref{BFS_on_H} shows the spanning forest found by BFS for a graph $G$.

The queue $Q$ for Figure \ref{BFS_on_H}  is 
\begin{center}
\begin{tabular}{|c|c|c|c|c|c|c|c|c|c|c|c|c|}
\hline
\hspace{.2cm} t\hspace{.2cm} &
\hspace{.2cm} 0\hspace{.2cm} &
\hspace{.2cm} 1\hspace{.2cm} &
\hspace{.2cm} 2\hspace{.2cm} &
\hspace{.2cm} 3\hspace{.2cm} &
\hspace{.2cm} 4\hspace{.2cm} &
\hspace{.2cm} 5\hspace{.2cm} &
\hspace{.2cm} 6\hspace{.2cm} &
\hspace{.2cm} 7\hspace{.2cm} &
\hspace{.2cm} 8\hspace{.2cm} &
\hspace{.2cm} 9\hspace{.2cm} &
\hspace{.2cm} 10\hspace{.2cm} &
\hspace{.2cm} 11\hspace{.2cm} \\ \hline 
Q & (1) & (3,4) & (4,8) & (8,7) & (7) & (6,9) & (9) & (2) & (5,10)  &
(10,11) & (11) & $\emptyset$ \\
\hline
\end{tabular}
\end{center}

For a spanning forest $F$ of $G$, let us say that an edge $e \in G-F$ is \emph{BFS-externally active} if $BFS(F\cup e)=F$. A crucial observation is 
made by Spencer \cite{Spencer:1997}: An edge $\{v,w\}$ can be added to $F$ without changing the spanning forest under the BFS if and only if the two 
vertices $v$ and $w$ have been present in the queue at the same time.  In our example of Figure \ref{BFS_on_H}, edges $\{3,4\}, \{4,8\}, \{7,8\}, 
\{6,9\}, \{5,10\}$ and $\{10, 11\}$ could be added back to $F$.  We write $\mathcal{E}(F)$ for the set of BFS-externally active edges.

\begin{prop}[Spencer] \label{eactive}
If $H$ is any subgraph and $F$ is any spanning forest of $G$ then 
$$
BFS(H)=F \text{ if and only if } F \subseteq H \subseteq F \cup \mathcal{E}(F).
$$
\end{prop} 

Now consider the Tutte polynomial. Note that if $BFS(H)=F$, then
 $c(H)=c(F)$. 
 So $\sigma(H)=c(F)-1$ and
$\sigma^*(H)=|E(H)|-|E(F)|=|\mathcal{E}(F) \cap H|$.  Hence if we fix
 a 
forest $F$ and sum over the corresponding
interval $[F, F \cup \mathcal{E}(F)]$, we have
$$
\sum_{H: BFS(H)=F}x^{\sigma(H)} y^{\sigma^*(H)} =x^{c(F)-1} \sum_{A \subseteq
\mathcal{E}(F)}y^{|A|} = x^{c(F)-1}(1+y)^{|\mathcal{E}(F)|}.
$$
Summing over all forests $F$, we get
\begin{eqnarray*}
t_G(1+x, 1+y)= \sum_{H \subseteq G }x^{\sigma(H)} y^{\sigma^*(H)}
=\sum_{F \subseteq G}  x^{c(F)-1}(1+y)^{|\mathcal{E}(F)|}. 
\end{eqnarray*} 
Or, equivalently, 
\begin{eqnarray} \label{Tutte_E(T)}
t_G(1+x,y)= \sum_{F \subseteq G}  x^{c(F)-1}y^{|\mathcal{E}(F)|}.
\end{eqnarray} 

To evaluate $\mathcal{E}(F)$, note that when applying BFS to a graph $H$, the queue $Q$ only depends on the spanning forest $F=BFS(H)$. Given a 
forest $F$, the processing order in $Q$ is a total order $<_Q=<_Q(F)$ on the vertices of $F$ satisfying the following condition: Let $T_1, T_2, 
\dots, T_k$ be the tree components of $F$ with minimal elements $r_1=1 < r_2 < \cdots < r_k$. Then (1) If $v$ is a vertex in tree $T_i$, $w$ is a 
vertex in tree $T_j$ and $i<j$, then $v<_Q w$. (2) Among vertices of each tree $T_i$, $r_i$ is minimal in the order $<_Q$. (3) For two non-root 
vertices $v,w$ in the same tree, $v<_Q w$ if $v^p <_Q w^p$. In the case $v^p=w^p$, $v<_Q w$ whenever $v<w$.

Comparing with the examples in \S3, we note that $<_Q$ is exactly the order $<_{bf,q}$ described in Example 5 of \S3, as \emph{breadth-first order 
with a queue}. Fix the choice function $\gamma=\gamma_{bf,q}$, the one associated to $<_{bf,q}$ and consider the maps $\Phi_{\gamma,G}$ and 
$\Psi_{\gamma, G}$. Given $F$, the condition that two vertices $v,w$ have been present at the queue $Q$ at the same time when applying BFS to $F$ is 
equivalent to $v^p <_{bf,q} w <_{bf,q} v$ or $w^p <_{bf,q} v <_{bf,q} w$.  That is, an edge is BFS-externally active if and only if it is an 
$F$-redundant edge of type 3, as defined in \S4.1.  It follows that $\mathcal{E}(F)=R_3(G;F)$.

Therefore by Prop. \ref{edgecategorization}, $$ |\mathcal{E}(F)| = |R_3(G;F)|= |E(G)|-|E(F)|-\bigg(\sum_{v:f(v)=-1} f(v)\bigg) -Rec(f), $$ where 
$f=\Psi_{\gamma, G}(F)$ is the corresponding \gmp.  Note that $|E(F)|=n-c(F)$, and $c(F)=r(f)$, where $r(f)$ is the number of roots of $f$. 
Therefore

\begin{thm} \label{new_Tutte} 
$$
t_G(1+x, y) = y^{|E(G)|-n} 
\sum_{f} x^{r(f)-1} y^{r(f)-Rec(f)-\left(\sum_{v: f(v) \neq \infty} 
f(v) \right)},
$$
where the sum is over all $G$-multiparking functions. 
\end{thm} 

For a \gmp \ $f$, where $G$ is a graph on $n$ vertices, we call the statistics $|E(G)|-n+r(f)-Rec(f)-\sum_{v: f(v) \neq \infty} f(v)$ the {\em 
reversed sum} of $f$, denote by $rsum(f)$. The name comes from the corresponding notation for classical parking functions, see, for example, 
\cite{KungYan:2003_2}.  Theorem \ref{new_Tutte} expresses Tutte polynomial in terms of generating functions of $r(f)$ and $rsum(f)$. In 
\cite{GesselSagan:1996} Gessel and Sagan gave a similar expression, in terms of $\mathcal{E}_{DFS}(F)$, the set of \emph{greatest-neighbor 
externally active} edges of $F$, which is defined by applying the greatest-neighbor depth-first search on subgraphs of $G$. Combining the result of 
\cite{GesselSagan:1996} (Formula 5), we have
\begin{eqnarray} 
xt_G(1+x, y) =\sum_{F \subseteq G} x^{c(F)}y^{|\mathcal{E}_{DFS}(F)|}
            =\sum_{F \subseteq G} x^{c(F)}y^{|\mathcal{E}(F)|}
            =\sum_{f \in \mpf} x^{r(f)}y^{rsum(f)}. 
\end{eqnarray}

That is, the three pairs of statistics, $(c(F),|\mathcal{E}_{DFS}(F)|)$ and $(c(F),|\mathcal{E}(F)|)$ for spanning forests, and $(r(f), rsum(f))$ 
for $G$-multiparking functions, are equally distributed.

\textsc{Remark}.  Alternatively, one can prove Theorem \ref{new_Tutte} by conducting \emph{neighbors-first search} (NFS), a tree traversal defined 
in \cite[\S6]{GesselSagan:1996}, and using $\gamma=\gamma_{df}$, the choice function associated with the depth-first search order, (c.f. Example 3, 
\S3).  Here the NFS is another algorithm that builds a spanning forest $F$ given an input graph $H$.  The following description is taken from 
\cite{GesselSagan:1996}.

\begin{itemize}
\item[NFS1] Let $F=\emptyset$. 
\item[NFS2] Let $v$ be the least unmarked vertex in $V$ and mark v. 
\item[NFS3] Search $v$ by marking all neighbors of $v$ that have not been 
marked and adding to $F$ all edges from $v$ to these vertices. 
\item[NFS4] Recursively search all the vertices marked in NFS3 in increasing 
order, stopping when every vertex that has been marked has also been searched.
\item[NFS5] If there are unmarked vertices, then return to
  NFS2. Otherwise, 
stop. 
\end{itemize} 

\begin{figure}[ht]
\centerline{\epsfxsize=5in \epsfbox{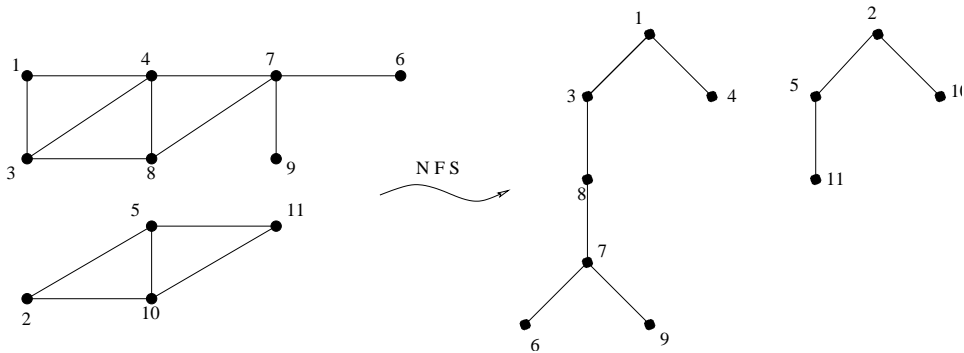}}
\caption{Spanning forest  found by NFS.}\label{NFS_on_H}
\end{figure}

The NFS searches vertices of $H$ in a depth-first manner but marks
children 
in a locally breadth-first manner. Figure \ref{NFS_on_H} shows 
the result of NFS, when applies to the graph on the left
of Figure \ref{BFS_on_H}.

Similarly, one defines $\mathcal{E}_{NFS}(F)$, the set of edges 
{\em externally active with respect to NFS}, to be
those edges $e \in G-F$ such that $NFS(F \cup e)=F$. 
Then Prop.  \ref{eactive} and Eq. \eqref{Tutte_E(T)} hold again
when we replace BFS with NFS, and $\mathcal{E}(F)$ with $\mathcal{E}_{NFS}(F)$.

Now let $\gamma=\gamma_{df}$ and use the bijections
 $\phi_{\gamma_{df}, G}$
 and $\Psi_{\gamma_{df}, G}$, one notices
again that an edge is externally active with respect to NFS if and
 only if 
it is $F$-redundant of type 3. And hence we
get another proof of Theorem \ref{new_Tutte}.

\vspace{.5cm} 

An interesting specialization of Theorem \ref{new_Tutte} is to
consider 
$t_G(1, y)$, the restriction to spanning trees
of $G$ and $G$-parking functions.  For a $G$-parking function $f$, 
or equivalently a $G$-multiparking function with
exactly one root (which is vertex 1), $r(f)=1$ and $Rec(f)=0$. 
Hence $rsum(f)=|E(G)|-n+1-\sum_{v \neq \infty} f(v)$. Thus
we obtain
$$
t_G(1,y)=\sum_{f \: : \: \text{$G$-parking functions} }y^{rsum(f)}. 
$$

An equivalent form of this result, in the language of sand-pile
models, 
was first proved by L\'opez \cite{Lopez:1997}
using a recursive characterization of Tutte polynomial. A bijective
proof 
was given by Cori and Le Borgne in
\cite{CoriLeborgne:2003} by constructing a one-to-one correspondence 
between trees with external activity $i$ (in
Tutte's sense) to recurrent configurations of level $i$, which is 
equivalent to $G$-parking functions with reversed
sum $i$. Our treatment here provides a new bijective proof.

In \cite{GesselSagan:1996} it is shown that, restricted to simple
graphs, 
the greatest-neighbor externally active
edges of $F$ are in one-to-one correspondence with certain inversions
of $F$. 
 For a simple graph $G$, view each tree
$T$ of $F$ as rooted at its smallest vertex.  An edge $\{u,v\}$ is 
greatest-neighbor externally active if and only if
$v$ is a descendant of $u$, and $w>v$ where $w$ is the child of $u$ 
on the unique $u-v$ path in $F$, (that is,
$u=w^p$). Call such a pair $\{w, v\}$ a $G$-inversion. 
And denoted by $Ginv(F)$ the number of $G$-inversions of the
forest $F$. Then we have the following corollary.

\begin{cor} \label{inv_rsum}
Let $\mathcal{F}_k(G)$ be the set of spanning forests of $G$ with exactly
$k$ tree components. And $\mathcal{MP}_k(G)$ be the set of $G$-multiparking
functions with $k$ roots. Then 
$$
\sum_{F \in \mathcal{F}_k(G)} y^{Ginv(F)} =\sum_{f \in \mathcal{MP}_k(G)} 
y^{rsum(f)}. 
$$
\end{cor} 

In particular, when $G$ is the complete graph $K_{n+1}$ and $k=1$, 
we have the well-known result on the
equal-distribution of inversions over labeled trees, and the reversed 
sum over all classical parking functions of
length $n$, (for example, see \cite{Kreweras,Stanley:1998})
$$
\sum_{T \text{ on } [n+1]} y^{inv(T)} 
=\sum_{\alpha \in P_n} y^{{n \choose 2}-\sum_{i=1}^n{\alpha_i}},
$$ 
where $P_n$ is the set of all (classical) parking functions of length $n$.

%
%

\section{Enumeration of \gmps \ and graphs}   


In this section we discuss some enumerative results on \gmps \ and
substructures of graphs.

\begin{thm}\label{number_gpf} 
The number of \gmps \ with $k$ roots equals the number of spanning 
forests of $G$ with $k$ components. In particular,
for connected graph $G$, the number of \gmps \ is $T_G(2,1)$. Among
them, 
those with an odd number of roots is counted
by $\frac 12 (T_G(2,1)+T_G(0,1))$, and those with an even number of
roots 
is counted by $\frac 12
(T_G(2,1)-T_G(0,1))$.
\end{thm} 

\begin{proof} 
The first two sentences follow directly from the bijections
constructed 
in \S2, and Theorem \ref{new_Tutte}. For the
third sentence, just note that $T_G(0,1)=\sum_{f} (-1)^{r(f)-1}$ is
the 
difference between the number of \gmps \ with
an odd number of roots, and those with an even number of roots.
\end{proof} 

Another consequence of Theorem \ref{new_Tutte} and its proof is an 
expression for the number of spanning subgraphs
with a fixed number of components and fixed number of edges, in terms 
of (BFS)-external activity and \gmps. It is a
generalization of the expectation formula in \cite{Spencer:1997},
which 
is the special case for complete graph
$K_{n}$.

\begin{thm} 
Let $G$ be a connected graph. The number $\gamma_{t,k}(G)$ of spanning 
subgraphs $H$ with $t$ components and
$V(G)-1+k$ edges is given by
$$
\gamma_{t,k}(G)=\sum_{F\in \mathcal{F}_t}{ \mathcal{E}(F)\choose k} =\sum_{f \in \mathcal{MP}_t(G)} {rsum(f) \choose k},
$$
where the first sum is over all spanning forests with $t$ components, 
and the second sum is over all \gmps \ with $t$
roots.
\end{thm} 

\begin{proof} 

For any spanning forest $F$ with $k$ components, the number of
spanning 
subgraphs $H$ with $V(G)-1+k$ edges such that
$BFS(H)=F$ is given by ${ \mathcal{E}(F)\choose k}$.
\end{proof} 

Next we give a new expression of the $t_{K_{n+1}}(x,y)$ in terms of classical parking functions. It enumerates the classical parking functions by 
the number of \emph{critical left-to-right maxima}. Given a classical parking function ${\bf b}=(b_1, \dots, b_n)$, we say that a term $b_i=j$ is 
\emph{critical} if in $\bf b$ there are exactly $j$ terms less than $j$, and exactly $n-1-j$ terms larger than $j$. For example, in ${\bf 
b}=(3,0,0,2)$, the terms $b_1=3$ and $b_4=2$ are critical. Among them, only $b_1=3$ is also a left-to-right maximum.

Let $\alpha(\bf b)$ be the 
number of critical left-to-right maxima in a classical parking
function $\bf b$. We have
\begin{thm} \label{opf} 
$$
t_{K_{n+1}}(x,y) =\sum_{{\bf b}\in P_n } x^{\alpha(\bf b)} y^{{n \choose
    2}-\sum_i b_i},
$$
where $P_n$ is the set of classical parking functions of length $n$.
\end{thm} 

\begin{proof} 
Let $F$ be a spanning forest on $[n+1]$ with tree components $T_1,
\dots, T_k$, where $T_i$ has minimal vertex $r_i$, and $r_1 <r_2 <
\cdots < r_k$. We define an operation $merge(F)$ which combines  
the trees $T_1, \dots, T_k$ by adding an 
edge between  $r_i$ with $w_{i-1}$ for each $i=2, ..., k$, 
where $w_{i-1}$ is the vertex of $T_{r-1}$ that is maximal under the 
order $<_{bf, q}$. Denote 
by $T_F=merge(F)$ the resulting tree. We observe
that for the forest $F$ and the tree $T_F$, the queue obtained by 
applying BFS are exactly the same. This implies that $F$ and $T_F$ 
have the same set of BFS-externally active edges. 

Conversely, given $T$ and an edge $e=\{w,v\} \in T$
where $w <_{bf, q} v$.  
We say the edge $e$ is \emph{critical} in $T$ if
$merge(T\setminus \{e\})=T$. 
Assume $T\setminus \{e\}=T_1 \cup T_2$ where $w \in T_1$ and $v \in
T_2$. 
By the  definition of the merge operation, $e$ is critical if and only
if $w$ is the maximal in $T_1$ under the order $<_{bf,q}$, and 
$v$ is vertex of  the lowest index in $T_2$. In terms of the queue
obtained by applying BFS to $T$, it is equivalent to
the following two conditions: (1)  There is a set $Q_i$ 
such that $Q_i=\{v\}$, and $v$ does not belong to any other $Q_i$.
(2) $v$ is of minimal index among the set of vertices processed 
after $v$.

Consider the maps $\Phi_{\gamma, G}$ and $\Psi_{\gamma, G}$
with $\gamma=\gamma_{bf, q}$ and $G=K_{n+1}$.
Let $f=\Psi_{\gamma, G}(T)$, and write  $f$ as a sequence $(f(2),f(3), \dots,
f(n+1))$. (There is no need to record $f(1)$, as $f(1)=\infty$
always.) Then
an edge $\{w, v\}$  is critical in $T$  if and only if 
(1) $f(v)$ is critical in the sequence $(f(2), \dots, f(n+1))$, 
and (2) $w > v$ for any vertex $w$ with $f(w)>f(v)$. 
That is, $f(v)$ is a left-to-right maximum in the sequence
 $(f(2),f(3), \dots, f(n+1))$. 
 
Now fix a spanning tree $T$ of $K_{n+1}$ and 
let $Merge(T)$ be the set of spanning forests $F$ such that 
$merge(F)=T$.   Then an $F \in Merge(T)$ can be obtained from $T$ by
removing any subset $A$ of critical edges, in which case 
$c(F)=c(T)+|A|=1+|A|$. This, combined with the fact that
 $\mathcal{E}(F)=\mathcal{E}(T)$, gives us
\begin{eqnarray} \label{sum-opf}
\sum_{F \in Merge(T)} x^{c(F)-1} y^{|\mathcal{E}(F)|}
=y^{|\mathcal{E}(T)|} \sum_A x^{|A|}, 
\end{eqnarray} 
where $A$ ranges over all subsets of critical edges of $T$. 
Under the correspondence $T \rightarrow f=\Psi_{\gamma, G}(T)$
and considering $f$ as a sequence $(f(2), \dots, f(n+1))$, 
$|\mathcal{E}(T)|$ is just ${n\choose 2}-\sum_{i=2}^{n+1} f(i)$, and
critical edges of $T$ correspond to critical left-to-right maxima of 
the sequence. 
Hence the sum in \eqref{sum-opf} equals 
$$
y^{|\mathcal{E}(T)|} (1+x)^{\alpha(f_T)}= (1+x)^{\alpha(f_T)} y^{{n
    \choose 2}-\sum_{i=2}^{n+1} f(i)}.
$$
Theorem \ref{opf} follows by summing over all trees on $[n+1]$. 
\end{proof}

Finally, we use the breadth-first search to re-derive the formula for the
number of subdigraphs of $G$, which was first
proved in \cite{GesselSagan:1996} using DFS, and extend the method to 
derive a formula for the number of subtraffics of $G$. 

Let $G$ be a graph. A {\em directed subgraph} or {\em subdigraph} of
$G$ 
is a digraph $D$ that contains up to one copy
of each orientation of every edge of $G$. Here for an edge $\{u,v\}$
of 
$G$ we permit both $(u,v)$ and $(v,u)$ to
appear in a subdigraph.

For any subdigraph $D$ of $G$, we apply the BFS to get a spanning
forest 
of $D$. The only difference from the subgraph
case is that when processing a vertex $x$, we only add those unvisited 
vertices $u$ such $(x,u)$ is an edge of $D$.

If digraph $D$ has BFS forest $F$, write $\mathcal{\vec F}^+(D)=F$. 
Note that we can view $F$ as an oriented spanning
forest, where each edge is pointing away from the root (i.e., the
minimal 
vertex) of the underlying tree component.  
Say a directed edge $\vec e \notin F$ is \emph{directed BFS externally 
active} with respect to $F$ if $\mathcal{\vec
F}^+(F\cup \vec e) =F$. Denote by $\mathcal{E}^+(F)$ the set of
directed
 BFS-externally active edges. Then we have the
following basic proposition, which is the analog in the undirected 
case.

\begin{prop} 
If $D$ is any subdigraph and $F$ is any spanning forest of $G$
then 
$$
\mathcal{\vec F}^+(D)=D \text{ if and only if } 
F \subseteq D \subseteq F \cup \mathcal{E}^+(F). 
$$
\end{prop} 

Now we characterize the directed BFS-externally active edges 
by the set $\mathcal{E}(F)$, the BFS-externally active edges 
for the undirected graph $G$. Let $\{u,v\}$ be an edge of $G$ with 
$u <_{bf, q} v$. If $\{u, v\} \in E(F)$, then the backward edge
$(v,u)$ can be added without changing the result of (directed) breadth-first
search, that is,  $(v,u) \in \mathcal{E}^+(F)$. 
If $\{u,v\}\in \mathcal{E}(F)$, then both 
$(u,v)$ and $(v,u)$ are in $\mathcal{E}^+(F)$. If $\{u,v\}$ is not 
in the forest $F$  or $\mathcal{E}(F)$, then $(v,u)$ is in 
$\mathcal{E}^+(F)$. Together we have
$$
|E(G)|=|\mathcal{E}^+(F)| -|\mathcal{E}(F)|. 
$$
Therefore 
\begin{thm} \label{subgraph} 
If $G$ has $n$ vertices, 
then 
\begin{eqnarray} 
\sum_D x^{c(D)}y^{|E(D)|} = xy^{n-1} (1+y)^{|E(G)|} \ t_G(1+\frac xy, 1+y),
\end{eqnarray} 
where the sum is over all subdigraphs of $G$.
\end{thm} 
\begin{proof} 
\begin{eqnarray*} 
 \sum_D x^{c(D)}y^{|E(D)|} & =& \sum_{F} \sum_{D: \mathcal{\vec F}^+(D)=F} 
 x^{c(D)}y^{|E(D)|} \\
 & = & \sum_F x^{c(F)} y^{|E(F)|} (1+y)^{|\mathcal{E}^+(F)|} \\
 & =&  y^n (1+y)^{|E(G)|} \sum_F \left( \frac xy\right)^{c(F)} (1+y)^{|\mathcal{E}(F)|} \\
 & = & xy^{n-1}  (1+y)^{|E(G)|} \ t_G(1+\frac xy, 1+y).
\end{eqnarray*}
\end{proof} 

Next we consider a  slightly complicated problem.
The {\em sub-traffic} $K$ of $G$, where $K$ is a
partially directed graph on $V(G)$, is obtained from $G$ by replacing 
each edge $\{u,v\}$ of $G$ by (a) $\emptyset$,
(b) a directed edge $(u,v)$, (c) a directed edge $(v,u)$, (d) two
directed
 edges $(u,v)$ and $(v,u)$, or (e) an
undirected edge $\{u,v\}$.  We proceed as we did before.  
For each subtraffic $K$, we apply the directed breadth-first
search to get a spanning forest $F$: The queue starts with the minimal 
vertex $1$. At each iteration, we take the
vertex $x$ at the head of the queue, remove $x$ from the queue, 
and add all unvisited vertices $u$ if $(x,u) \in E(K)$
or $\{x, u\} \in E(K)$. Add the directed edge $(x,u)$ to the forest
$F$ 
if $(x,u) \in E(K)$. Otherwise, add the
undirected edge $\{x,u\}$ to $F$. The output is a forest $[n]$
 in which each edge is either a directed edge oriented
away from the minimal vertex of the underlying tree, or an undirected
edge. 
Let $A$ be the set of directed edges.
Denote by $(F, A)$ the output forest and write $BFS(K)=(F, A)$.  
Note that $(F,A)$ is itself a sub-traffic of $G$.

Given a pair $(F,A)$ with directed edges $A \subseteq E(F)$,
we have the following characterization of edges  that can be added 
to $(F, A)$, without changing the BFS result, (i.e., $BFS( (F,A) \cup
e) =BFS(F,A)$.)  
\begin{enumerate} 
\item  For each directed edge $(u,v)$ in $A\subseteq E(F)$, we can add back
$(v,u)$ without changing the result of the spanning forest. 
\item  For each BF-externally active edge $\{u,v\}$ of $F$, we can add
  back any one of $(u,v), (v, u)$ and  $\{u,v\}$, or both $(u,v)$ and
  $(v,u)$  at the same time. 
\item  For each edge not in $F \cup \mathcal{E}(F)$, we can add back 
one of the undirected edge $\{u,v\}$ and the direct  $(u,v)$ if $u$  
is processed after $v$ in the queue. 
\end{enumerate} 

There is no further restriction on how the edges can be added back
 in addition to the above mentioned cases. Then we have
\begin{thm} 
Let $G$ be a connected graph. Then
\begin{eqnarray}~\label{newtuttesum}
\sum_K x^{c(K)} y^{|E(K)|} = x (y^2+2y)^{n-1} (1+2y)^{|E(G)|-n+1} \ t_G
(1+\frac{x(1+2y)}{y(2+y)}, \frac{1+3y+y^2}{1+2y}),
\end{eqnarray} 
where the sum is over all subtraffic of $G$. 
\end{thm} 
\begin{proof} 
$$
\sum_K x^{c(K)} y^{|E(K)|} =\sum_F \sum_{A \subseteq E(F)} \sum_{K:
  BFS(K)=(F,A)} x^{c(K)} y^{|E(K)|},
$$
where $F$ is over all spanning forests of $G$, and $A$ is a subset
of the edges of $F$. 
A subtraffic $K$ has $BFS(K)=(F, A)$ if and only if it is obtained from 
$F$ by adjoining some edges as described in the preceding three 
cases. Considering the contribution of each type, we have
\begin{eqnarray*}
 \sum_{K:  BFS(K)=(F,A)} x^{c(K)} y^{|E(K)|}&  =& x^{c(F)} y^{|E(F)|} 
(1+y)^{|A|} (1+3y+y^2)^{|\mathcal{E}(F)|}
 (1+2y)^{|E(G)|-|E(F)|-|\mathcal{E}(F)|} \\
& = &  x^{c(F)} \left(\frac{y}{1+2y}\right)^{|E(F)|} (1+2y)^{|E(G)|}
 (1+y)^{|A|} \left( 
\frac{1+3y+y^2}{1+2y}\right) ^{|\mathcal{E}(F)|}.
\end{eqnarray*} 
Hence
\begin{eqnarray*}
 \sum_K x^{c(K)} y^{|E(K)|} 
&=&\sum_F  x^{c(F)} \left(\frac{y}{1+2y}\right)^{|E(F)|} (1+2y)^{|E(G)|}
 \left( \frac{1+3y+y^2}{1+2y}\right) ^{|\mathcal{E}(F)|} \sum_{A
 \subseteq E(F)} (1+y)^{|A|} \\
& =& \sum_F  x^{c(F)} \left(\frac{y}{1+2y}\right)^{|E(F)|} (1+2y)^{|E(G)|}
 \left( \frac{1+3y+y^2}{1+2y}\right) ^{|\mathcal{E}(F)|}
 (2+y)^{|E(F)|} \\
& = & \left( \frac{y(2+y)}{1+2y}\right)^n (1+2y)^{|E(G)|} \sum_F \left(
 \frac{x(1+2y)}{y(2+y)}\right)^{c(F)}  \left(
 \frac{1+3y+y^2}{1+2y}\right) ^{|\mathcal{E}(F)|} \\
& = & x (y^2+2y)^{n-1} (1+2y)^{|E(G)|-n+1} \
 t_G\left(1+\frac{x(1+2y)}{y(2+y)}, \frac{1+3y+y^2}{1+2y}\right).
\end{eqnarray*} 
 \end{proof} 

By evaluating equation \ref{newtuttesum} at $x=y=1$, we derive a new evaluation of the Tutte polynomial that counts the number of subtraffics $K$ on 
$G$.

\begin{cor} The number of subtraffics on $G$ is equal to $3^{|E(G)|} t_G(2, \frac{5}{3})$.

\end{cor}

\section*{Acknowledgments}

The authors thank Robert Ellis and Jeremy Martin for helpful
discussions 
and comments. We also thank Ira Gessel for
helpful comments on Tutte polynomials, and for sharing the unpublished 
portion of a preprint of
\cite{GesselSagan:1996} with us.


\end{document}